\documentclass[10pt]{amsart}
\usepackage{amscd,amssymb,amsmath, array,latexsym,amsbsy}

\allowdisplaybreaks[2]
\begin{document}
%%%%%%%%%%%%%%%%%%%%%%%%%%%%%%%%%%

\newcommand{\go}{G^0}
\newcommand{\cc}{C_{c}}
\newcommand{\thmref}[1]{Theorem~\ref{#1}}
\newcommand{\defref}[1]{Definition~\ref{#1}}
\newcommand{\lemref}[1]{Lemma~\ref{#1}}
\newcommand{\propref}[1]{Prop\-o\-si\-tion~\ref{#1}}
\newcommand{\corref}[1]{Cor\-ol\-lary~\ref{#1}}

%%%%%%%%%%%%%%%%%%%%%%%%%%%%%%%%%%
\newcommand{\cs}{\ensuremath{C^{*}}}
\def\supp{\operatorname{supp}}
\def\tr{\operatorname{tr}}
\def\lt{\operatorname{{lt}}}
\def\interior{\operatorname{Int}} % interior
\def\Ind{\operatorname{Ind}} % induced representation
\def\Prim{\operatorname{Prim}}
%%%%%%%%%%%%%%%%%%%%%%%%%%%%
\def\H{\mathcal{H}} %Hilbert space
\def\K{\mathcal{K}} % Compacts

\def\C{\mathbb{C}}
\def\T{\mathbb{T}}
\def\Z{\mathbb{Z}}
\def\N{\mathbb{N}}
\def\R{\mathbb{R}}
\def\P{\N\setminus\{0\}}
%%%%%%%%%%%%%%%%%%%%%%%%%%%%%%%%%%
\newtheorem{thm}{Theorem}[section]
\newtheorem{cor}[thm]{Corollary}
\newtheorem{prop}[thm]{Proposition}
\newtheorem{lemma}[thm]{Lemma}
\theoremstyle{definition}
\newtheorem{defn}[thm]{Definition}
\newtheorem{remark}[thm]{Remark}
\newtheorem{remarks}[thm]{Remarks}
\newtheorem{claim}[thm]{Claim}
\newtheorem{example}[thm]{Example}
\numberwithin{equation}{section}
%%%%%%%%%%%%%%%%%%%%%%%%%%%%%%%%%%%
\title
[\boldmath Principal groupoid $C^*$-algebras with bounded trace]
{Principal groupoid $C^*$-algebras with bounded trace}

\author[Clark]{Lisa Orloff Clark}
\address{Department of Mathematical Sciences
\\Susquehanna University
\\Selinsgrove,  PA 17870
\\USA
} \email{clarklisa@susqu.edu}

\author[an Huef]{Astrid an Huef}
\address{School of Mathematics and Statistics
\\University of New South Wales
\\Sydney, NSW 2052
\\Australia}
\email{astrid@unsw.edu.au}

\date{August 22, 2006}

\begin{abstract}
Suppose $G$ is a second countable, locally compact, Hausdorff,
principal groupoid with a fixed left Haar system.  We define
a notion of integrability for groupoids, and show $G$ is integrable if and only
if the groupoid $C^*$-algebra $\cs(G)$ has bounded trace.
\end{abstract}

\thanks{This research was supported by the Australian Research Council and an AWM-NSF Mentoring Travel Grant.}

\maketitle
\section{Introduction}
Let $H$ be a locally compact, Hausdorff group acting continuously on a locally compact, Hausdorff space $X$, so that $(H,X)$ is a transformation group.  A lovely theorem of Green says that if $H$ acts freely on $X$, then the associated transformation-group $C^*$-algebra $C_0(X)\rtimes H$ has continuous trace if and only if the action of $H$ on $X$ is proper \cite[Theorem~17]{green}.
Muhly and Williams defined a notion of proper groupoid, and proved that for principal groupoids $G$, the groupoid $C^*$-algebra $C^*(G)$ has continuous trace if and only if the groupoid is proper \cite[Theorem~2.3]{MW90}.  Of course, when $G=H\times X$ is the transformation-group groupoid, then $G$ is proper if and only if the action of $H$ on $X$ is proper.

In \cite{rieffel} Rieffel introduced a notion of an integrable action of a group $H$ on a $C^*$-algebra $A$. This notion of integrability for $A=C_0(X)$ turned out to characterize when $C_0(X)\rtimes H$, arising from a free action of $H$ on $X$, has bounded trace \cite[Theorem~4.8]{aH02}. In this paper we define a notion of integrability for groupoids (see Definition~\ref{defn-int}) which, when $G=H\times X$ is the transformation-group groupoid, reduces to an integrable action of $H$ on $X$ (see Example~\ref{ex-transf}).  We then prove that for principal groupoids $G$, $C^*(G)$ has bounded trace if and only if the groupoid is integrable (see Theorem~\ref{thm-main}).  This theorem is thus very much in the spirit of \cite[Theorem~2.4]{MW90}, \cite[Theorem~7.9]{C}, \cite[Theorem~4.1]{C} (see also \cite[Corollary~5.9]{AS}) and \cite[Theorem~5.3]{C}, which characterize when principal-groupoid $C^*$-algebras are, respectively continuous-trace, Fell, CCR and GCR $C^*$-algebras.
The key technical tools used to prove Theorem~\ref{thm-main} are first, a homeomorphism of the spectrum of $C^*(G)$ onto the orbit space \cite[Proposition~5.1]{C}, and second, a generalisation to groupoids of the notion of $k$-times convergence in the orbit space of a transformation group from \cite{AD}.

\section{Preliminaries}

Let $A$ be a $C^*$-algebra. An element $a$ of the
positive cone $A^+$ of $A$ is called a bounded-trace element if the
map $\pi \mapsto \tr(\pi(a))$ is bounded on the spectrum $\hat A$ of
$A$; the linear span of the bounded-trace elements is a  two-sided $*$-ideal in $A$.  We say $A$  
has \emph{bounded trace} if the ideal  of (the span of) the
bounded-trace elements is dense in A. 

Throughout, $G$ is a   locally compact, Hausdorff
groupoid; in our main results $G$ is assumed to be second-countable and principal.
We denote the unit space of $G$  by $G^0$, and the range and
source maps $r,s:G\to G^0$ are $r(\gamma)=\gamma\gamma^{-1}$ and
$s(\gamma)=\gamma^{-1}\gamma$, respectively. We let $\pi:G\to
G^0\times G^0$ be the map $\pi(\gamma)=(r(\gamma),s(\gamma))$; recall that $G$ is  principal if $\pi$
is injective. In order to define the groupoid $\cs$-algebra, we also
assume that $G$ is equipped with a fixed left Haar system:
a set  $\{\lambda^x:x\in G^0\}$ of non-negative Radon measures on
$G$ such that
\begin{enumerate}
\item $\supp\lambda^x=r^{-1}(\{x\})$;
\item for $f\in C_c(G)$, the function $x\mapsto \int f\, d\lambda^x$ on $G^0$ is in $C_c(G^0)$; and
\item for $f \in \cc(G)$ and $\gamma\in G$, the following equation holds:
\begin{equation*} \int f(\gamma \alpha) \,d\lambda^{s(\gamma)}(\alpha)=\int f(\alpha)\,d
 \lambda^{r(\gamma)}(\alpha).
 \end{equation*}
\end{enumerate}
Condition (3) implies that $\lambda^{s(\gamma)}(\gamma^{-1}
E)=\lambda^{r(\gamma)}(E)$ for measureable sets $E$.  The collection
$\{\lambda_x:x\in G^0\}$, where $\lambda_x(E):=\lambda^x(E^{-1})$, gives a fixed right Haar system such that
the measures are supported on $s^{-1}(\{x\})$ and
\[
\int f(\gamma \alpha) \ d\lambda_{r(\alpha)} = \int f(\gamma) \
d\lambda_{s(\alpha)}
\]
for $f\in C_c(G)$ and $\gamma\in G$.   We will move freely between
these two Haar systems.

If $N\subseteq G^0$ then the \emph{saturation} of $N$ is
$r(s^{-1}(N))=s(r^{-1}(N))$. In particular, we call the saturation
of $\{x\}$ the \emph{orbit} of $x\in G^0$ and denote it by $[x]$.

If $G$ is principal and all the orbits are locally closed, then
by \cite[Proposition~5.1]{C} the orbit space $\go/G=\{[x]:x\in G^0\}$  and the
spectrum $\cs(G)^\wedge$  of the groupoid $C^*$-algebra $C^*(G)$ are homeomorphic.
This homeomorphism is induced by the the map $x\mapsto L^x:G^0\to C^*(G)^\wedge$,
where $L^x: C^*(G)\to B(L^2(G,\lambda_x))$  is given by
\[L^x(f)\xi(\gamma) = \int f(\gamma
\alpha)\xi(\alpha^{-1})d\lambda^x(\alpha)
\]
for $f\in C_c(G)$ and $\xi\in L^2(G,\lambda_x)$.

\section{Integrable Groupoids and convergence in the orbit space of a groupoid}

The following definition is motivated by the notion of  integrable action of a locally compact, Hausdorff group on a space from \cite[Definition~3.2]{aH02}.
\begin{defn}\label{defn-int}
A locally compact, Hausdorff groupoid $G$ is \emph{integrable} if for every compact subset $N$ of $G^0$,
\begin{equation}\label{eq-defn-integrable}
\sup_{x\in N}\{\lambda^x(s^{-1}( N))\}<\infty,
\end{equation}
or, equivalently, $\sup_{x\in N}\{\lambda_x(r^{-1}(N))\}<\infty$.
\end{defn}

\begin{remarks}\label{rem-integr}
 %Since the support of the measure $\lambda^x$ is $r^{-1}(\{x\})$,
%\[\lambda^x(s^{-1}(N))=\lambda^x(s^{-1}(N)\cap r^{-1}(\{x\})).\]

 We could have taken the supremum in   \eqref{eq-defn-integrable}
over the whole unit space, that is \[\sup_{x\in G^0}\{\lambda^x(s^{-1}( N))\}=
\sup_{x\in N}\{\lambda^x(s^{-1}( N))\}.\]
 To see this,  first note that if $y$ is not in the saturation $r(s^{-1}(N))=s(r^{-1}(N))$
 of $N$ then $s^{-1}(N)\cap r^{-1}(\{y\})=\emptyset$, and hence $\lambda^y(s^{-1}(N))=0$.  Second, if $y$ is in the saturation
 of $N$ then there exists a $\gamma\in G$ such that $s(\gamma)=y$ and $r(\gamma)\in N$.  Then
\[r^{-1}(\{y\})\cap s^{-1}(N)=\gamma^{-1}\gamma\big( r^{-1}(\{y\})\cap s^{-1}(N)\big)=
\gamma^{-1}\big(r^{-1}(\{r(\gamma)\})\cap s^{-1}(N)\big),
\]
and now
\begin{align*}
\lambda^y(s^{-1}(N))&=\lambda^y\big( r^{-1}(\{y\})\cap s^{-1}(N) \big)
=
\lambda^{r(\gamma^{-1})}\big( r^{-1}(\{r(\gamma)\})\cap s^{-1}(N) \big)\\
&=
\lambda^{s(\gamma^{-1})}\big(  r^{-1}(\{r(\gamma)\})\cap s^{-1}(N)   \big)=\lambda^{r(\gamma)}(s^{-1}(N))
\end{align*}
with $r(\gamma)\in N$.
\end{remarks}

\begin{example}\label{ex-transf}
Let $(H,X)$ be a  locally compact, Hausdorff transformation group
with $H$ acting on the left of the space $X$. Then $G=H\times X$ with
\[
G^2=\{((h,x),(k,y))\in G\times G \colon y=h^{-1}\cdot x\}
\]
and operations $(h,x)(k, h^{-1}\cdot x)=(hk,x)$ and
$(h,x)^{-1}=(h^{-1}, h^{-1}\cdot x)$ is called the
\emph{transformation-group groupoid}.  We identify the unit space
$\{e\} \times X$ with $X$, and then the range and source  maps $r,s:G\to X$
are  $s(h,x)=h^{-1}\cdot x$ and $r(h,x)=x$. If $\delta_x$ is the
point-mass measure on $X$ and $\mu$ is a left Haar measure on $H$,
then $\{\lambda^x:=\mu\times\delta_x\colon x\in X\}$ is a left Haar
system for $G$. Now
\[
\lambda^x(s^{-1}( N))=\mu(\{h\in H:h^{-1}\cdot x\in N\})
\]
and hence
\[
\sup_{x\in N}\{\lambda^x(s^{-1}( N))\}=\sup_{x\in N}\{\mu(\{h\in H:h^{-1}\cdot x\in N\})\},
\]
that is,  Definition~\ref{defn-int} reduces to \cite[Definition~3.2]{aH02}.
\end{example}

The following characterization of integrability will be important
later. In the case of a transformation-group groupoid,
Lemma~\ref{lem-integrable} reduces to a special case of
\cite[Lemma~3.5]{AD}.

\begin{lemma}\label{lem-integrable}
Let $G$ be a  locally compact, Hausdorff groupoid. Then $G$ is
integrable if and only if,
for each $z\in G^0$, there exists an open neighborhood of $U$ of $z$
in $G^0$ such that
\[
\sup_{x\in U}\{\lambda^x(s^{-1}(U))\}< \infty.
\]
\end{lemma}
\begin{proof}
The proof is exactly the same as the proof of \cite[Lemma~3.5]{AD}.
\end{proof}
If a groupoid fails to be integrable, there exists a $z\in
G^0$ such that
\[
\sup_{x\in U}\{\lambda^x(s^{-1}(U))\}=\infty
\]
for every open neighborhood  $U$ of $z$; we then say that the
\emph{groupoid fails to be integrable at $z$}.

It is evident from \cite{AD, AaH} that integrability and $k$-times
convergence in the orbit space of a transformation group are closely
related. Moreover, Lemma~2.6 of \cite{MW90} says that, if a
principal groupoid fails to be proper and the orbit space $G^0/G$ is
Hausdorff, then there exists a sequence that converges $2$-times in
$G^0/G$ in the sense of Definition~\ref{defn-ktimes}.

\begin{defn}\label{defn-ktimes}
A sequence $\{x_n\}$ in the unit space of a groupoid $G$ converges $k$-times in $G^0/G$ to $z\in G^0$ if there exist $k$ sequences
\[
\{ \gamma_n^{(1)}\},\  \{\gamma_n^{(2)}\},\dots,\{\gamma_n^{(k)}\}\subseteq G
\]
such that
\begin{enumerate}
\item $r(\gamma_n^{(i)})\to z$ as $n\to\infty$ for $1\leq i\leq k$;
\item  $s(\gamma_n^{(i)})=x_n$ for $1\leq i\leq k$;
\item if $1\leq i<j\leq k$ then $\gamma_n^{(j)}(\gamma_n^{(i)})^{-1}\to\infty$ as $n\to\infty$, in the sense that
 $\{\gamma_n^{(j)}(\gamma_n^{(i)})^{-1}\}$  admits no convergent subsequence.
\end{enumerate}
\end{defn}

\begin{remarks}
(a) Condition (2) in Definition~\ref{defn-ktimes} is needed so that the composition in (3) makes sense.

(b) Definition~\ref{defn-ktimes} does not require that $x_n\to z$,
but as in the transformation-group case
(\cite[Definition~2.2]{AaH}), this can be arranged by changing the
sequence which converges $k$-times: replace $x_n$ by
$r(\gamma_n^{(1)})$ and replace $\gamma_n^{(j)}$ by
$\gamma_n^{(j)}(\gamma_n^{(1)})^{-1}$.

(c) Part (3) of Definition~\ref{defn-ktimes} means
$\gamma_n^{(j)}(\gamma_n^{(i)})^{-1}$ is eventually outside every compact set.
In particular, if $LL^{-1}$ is compact, $L\gamma_n^{(i)}\cap L\gamma_n^{(j)}=\emptyset$
 eventually.
\end{remarks}

\begin{example}\label{ex-2times}
Let $G=H\times X$ be a transformation-group groupoid (see Example~\ref{ex-transf}) and suppose that $\{x_n\}\subseteq G^0$ is a sequence converging $2$-times in $G^0/G$ to $z\in G^0$.  Then there exist two sequences
\[\{\gamma_n^{(1)}\}=\{(s_n,y_n)\}\quad\text{and}\quad\{\gamma_n^{(2)}\}=\{(t_n,z_n)\}\]
in $G$ such that
(1)  $y_n\to z$ and $z_n\to z$;
(2) $s_n^{-1}\cdot y_n=x_n$ and $t_n^{-1}\cdot z_n=x_n$; and
(3) $(t_ns_n^{-1},z_n)\to\infty$ as $n\to\infty$.
To see that the sequence $\{x_n\}$ converges $2$-times in $X/H$ to $z$ in the sense of \cite[\S4]{AaH},
consider the two sequences $\{s_n\}$ and $\{t_n\}$ in $H$. We have  $s_n\cdot x_n\to z$ and $t_n\cdot x_n \to z$ using (1) and (2). Also, since $z_n\to z$ by (1), (3) implies that $t_ns_n^{-1}\to\infty$ in $H$.
\end{example}

In \S\ref{sec-four} we will prove that a principal groupoid $G$ is integrable if and only if $C^*(G)$ has bounded trace.  For the ``only if'' direction we will need to know that the orbits are locally closed so that \cite[Proposition~5.1]{C} applies and  $x\mapsto L^x$ induces a homeomorphism of $G^0/G$ onto $C^*(G)^\wedge$; Corollary~\ref{cor-orbits} below establishes that if $G$ is integrable then the orbits are locally closed. We will prove the contrapositive of the ``if'' direction, and a key observation for the proof is  Proposition~\ref{prop-int-ktimes}: if a groupoid fails to be integrable at some $z$ then there is a non-trivial sequence $\{x_n\}$ which converges $k$-times in $G^0/G$ to $z$, for every $k\in \P$. The next three lemmas are needed to establish Corollary~\ref{cor-orbits} and Proposition~\ref{prop-int-ktimes}.

Recall that a neighborhood $W$ of $\go$ is called \emph{conditionally compact} if the
sets $WV$ and $VW$ are relatively compact for every compact set $V$
in $G$.

\begin{lemma}\label{lem-measure}
Let $G$ be a  second countable, locally compact, Hausdorff groupoid.
\begin{enumerate}
\item   Let $z\in G^0$ and let $K$ be a relatively  compact neighborhood of $z$ in $G$.  There exist
$a\in\R$ and a neighborhood $U$ of $z$ in $G^0$ such that $0<a\leq
\lambda_x(K)$ for all $x\in U$.

\item Let $Q$ be a conditionally compact neighborhood in $G$. Given any relatively compact neighborhood
$V$  in $G^0$ such that $QV\neq\emptyset$, there exists $c\in\R$  such that $c>0$ and  $\lambda_x(Q)\leq
c$ for all $x\in V$.
\end{enumerate}
\end{lemma}
\begin{proof} (1) Suppose not. Let $\{U_i\}$ be a decreasing sequence of open 
neighborhoods of $z$ in $G^0$.  There exists an increasing seqence $i_1<i_2<\cdots<i_n<\cdots$ and 
$x_n\in U_{i_n}$ such that $\lambda_{x_n}(K)<1/n$ for each $n\geq 1$. Note that $x_n\to z$.

Let $f\in C_c(G)$
such that $0\leq f\leq 1$, $f(z)=1$ and $\supp f\subseteq K$; note
that $\int  f(\gamma) \, d\lambda_z(\gamma)>0$. By the continuity
of the Haar system,
\[
\frac{1}{n}>\lambda_{x_n}(K)\geq\int f(\gamma) \, d\lambda_{x_n}(\gamma)\to \int  f(\gamma) \,
d\lambda_z(\gamma) \text{\ as $n\to\infty$,}
\]
which is impossible since the left-hand side converges to $0$ and $\int f(\gamma) \, d\lambda_z(\gamma)>0$.
\medskip

(2) Let $V$ be any relatively compact neighborhood  in $G^0$ such that $QV\neq\emptyset$.
Let $f\in C_c(G)$ such that $0\leq f\leq 1$ and $f$ is identically
one on the relatively compact subset $QV$.  The function
$w\mapsto\int  f(\gamma)\, d\lambda_w(\gamma)$ is in $C_c(G^0)$, so
achieves a maximum $c>0$.  Then, for $x\in V$,
\[
\lambda_x(Q)=\lambda_x(Qx)\leq\int  f(\gamma)\,
d\lambda_x(\gamma)\leq c.
\]
\end{proof}

\begin{lemma} \label{lem-z}
Let $G$ be a  locally compact,
Hausdorff groupoid and let $z \in \go$.
Suppose the orbit $[z]$ is not locally closed and let $V$
be an open neighborhood of $z$ in $G^0$. For every compact neighborhood $K$
of $z$ in $G$ there exists $\gamma_K\in G\setminus K$ such that
$s(\gamma_K)=z$ and $r(\gamma_K)\in V$.
\end{lemma}

\begin{proof}
Suppose not, that is, suppose that there exist an open neighborhood $V$ of $z$ in $G^0$ and a compact neighborhood $K$ of $z$ in $G$ such that
for all $\gamma\notin K$, either $s(\gamma)\neq z$ or $r(\gamma)\notin V$.
% If $s(\gamma)\neq z$ for all $\gamma\notin K$ then $s^{-1}(\{z\})$
% is a closed subset of the compact set $K$ and hence is compact. But
% now $[z]=r(s^{-1}(\{z\}))$ is also compact and hence closed, a
% contradiction. So $s^{-1}(\{z\})\cap (G\setminus K)\neq \emptyset$.
Since the orbit $[z]$ is not locally closed, $(\overline{[z]}\setminus[z])\cap V\neq\emptyset$. Let $y\in (\overline{[z]}\setminus[z])\cap V$. Then there exists
$\{\gamma_i\}\subseteq G$ such that $s(\gamma_i)=z$, $r(\gamma_i)\to
y$.  Moreover, since $y\in V$ and $V$ is open,  $r(\gamma_i)\in V$ eventually.  So by assumption,
$\gamma_i\in K$ eventually. By passing to a subsequence and
relabeling we may assume that $\gamma_i\to\gamma\in K$.    But now
$s(\gamma)=z$ and $r(\gamma)=y$, which means that $y\in[z]$, a
contradiction.
\end{proof}

Lemma~\ref{lem-orbits} generalises \cite[Lemma~2.1]{AaH} from transformation groups to principal groupoids.

\begin{lemma}\label{lem-orbits}
Let $G$ be a second countable, locally compact, Hausdorff, principal groupoid and let
$z\in G^0$. Then the following four conditions are equivalent.

\begin{enumerate}
\item the orbit $[z]$ is not locally closed;
\item for every $k\in\P$, the sequence $\{z,z,\dots\}$ converges $k$-times in $G^0/G$ to $z$;
\item for every open neighborhood $V$ of $z$ in $\go$, $\lambda_z(r^{-1}(V))=\infty$;
\item for every open neighborhood $V$ of $z$ in $\go$,
$r^{-1}(V)\cap s^{-1}(\{z\})$ is not relatively
compact.
\end{enumerate}
\end{lemma}
\begin{proof}
(1) $\Longrightarrow$ (2). Suppose that $[z]$ is not open in
$\overline{[z]}$ and fix $k\in\P$. Then, for every neighborhood $W$ of
$z$,
\[
(\overline{[z]}\setminus[z])\cap W\neq\emptyset.
\]
Let $\{V_n\}$ be a decreasing sequence of basic open neighborhoods
of $z$ in $G^0$, and $K_n$ an increasing sequence of compact
neighborhoods of $z$ in $G$ such that $G=\cup_{n\geq
1}\interior(K_n)$.

Now we (temporarily) fix $n$ and choose $\gamma_n^{(1)},\dots,
\gamma_n^{(k)}$ as follows.  Set $\gamma_n^{(1)}=z$. By
Lemma~\ref{lem-z} there exists $\gamma_n^{(2)}\in G\setminus
(K_n\gamma_n^{(1)})$ such that $s(\gamma_n^{(2)})=z$ and
$r(\gamma_n^{(2)})\in V_n$. Next, note that $K_n\gamma_n^{(1)}\cup
K_n\gamma_n^{(2)}$ is compact, so by Lemma~\ref{lem-z} there exists
\[
\gamma_n^{(3)}\in G\setminus (K_n\gamma_n^{(1)}\cup K_n\gamma_n^{(2)})
\]
with $s(\gamma_n^{(3)})=z$ and $r(\gamma_n^{(3)})\in V_n$.  Continue.

By construction, $r(\gamma_n^{(i)})\to z$ for $1\leq i\leq k$ as
$n\to\infty$. Also, if $1\leq i<j\leq k$ then $\gamma_n^{(j)}\notin \cup_{i=1}^{j-1} K_n\gamma_n^{(i)}$; in particular
\[
\gamma_n^{(j)}(\gamma_n^{(i)})^{-1}\notin K_n.
\]
Now let $Q$ be any compact set in $G$.  Choose $n_0$ such that
$Q\subseteq K_n$ whenever $n>n_0$. Since
$\gamma_n^{(j)}(\gamma_n^{(j)})^{-1}\notin K_n$ we have
$\gamma_n^{(j)}(\gamma_n^{(j)})^{-1}\notin Q$ whenever $n>n_0$.
Hence $\gamma_n^{(j)}(\gamma_n^{(j)})^{-1}\to\infty$ as
$n\to\infty$. Thus $\{z,z,\dots\}$ converges $k$-times in
$G^{(0)}/G$ to $z$. Finally, note that $k$ was arbitrary.

\medskip
(2) $\Longrightarrow$ (3). Suppose that for every $k\in\P$, the
sequence $\{z,z,\dots\}$ converges $k$-times in $G^0/G$ to $z$. Let
$V$ be an open neighborhood of $z$ in $G^{(0)}$ and let $M\in \P$.
There exists an open neighborhood $U$ of $z$ in $G^{0}$
 and a compact symmetric neighborhood $K$ of $z$ in $G$ such that $r(KU)\subseteq V$.
 By Lemma~\ref{lem-measure}(1) and by shrinking $U$ if necessary, there exists $c>0$ such
 that $\lambda_u(K)\geq c$ for $u\in U$. Choose $k\in\P$ such that $kc>M$. Next, choose
 $k$ sequences $\{\gamma_n^{(i)}\}\subseteq G$ such that, for
\begin{enumerate}
\item[(a)] $s(\gamma_n^{(i)})=z$ and $r(\gamma_n^{(i)})\to z$ as $n\to\infty$ for $1\leq i\leq k$, and
\item[(b)] if $1\leq i<j\leq k$ then $\gamma_n^{(j)}(\gamma_n^{(i)})^{-1}\to\infty$ as $n\to\infty$.
\end{enumerate}
So there exists $n_0$ such that
$r(\gamma_{n_0}^{(i)})\in U$ for $1\leq i\leq k$ and $\gamma_{n_0}^{(j)}(\gamma_{n_0}^{(i)})^{-1}\in G\setminus K^2$.
Now  $r(K\gamma_{n_0}^{(i)})\subseteq r(KU)\subseteq V$ and $K\gamma_{n_0}^{(i)}\cap K\gamma_{n_0}^{(j)}=\emptyset$
when $i\neq j$.
By our choice of $n_0$ we have
\[
\lambda_z(r^{-1}(V))
%&=
%\lambda_z(\{\gamma: s(\gamma)=z \text{\ and\ } r(\gamma)\in V\})\\
\geq\lambda_z\Big(\bigcup_{i=1}^k K\gamma_{n_0}^{(i)}\Big)=\sum_{i=1}^k \lambda_z(K\gamma_{n_0}^{(i)})
=\sum_{i=1}^k \lambda_{r(\gamma_{n_0}^{(i)})}(K)\geq kc>M.
\]
Since $M$ was arbitrary, $\lambda_z(r^{-1}(V))=\infty$.

\medskip
(3) $\Longrightarrow$ (4) If $\lambda_z(r^{-1}(V)) =\infty$ then
\[
\lambda_z(r^{-1}(V)\cap s^{-1}(\{z\}))  =\infty
\]
because the support of $\lambda_z$ is $s^{-1}(\{z\})$.
Since $\lambda_z$ is a Radon measure $r^{-1}(V)\cap s^{-1}(\{z\})$ cannot
be relatively compact.

\medskip
(4) $\Longrightarrow$ (1)  We prove the contrapositive.  Assume
$[z]$ is locally closed in $G$.  Then there exists an open set $U$
in $G$ such that $U \cap \overline{[z]}=[z]$. Let $N$ be a compact
neighborhood of $z$ in $G$ such that $N \subset U$ (such $N$ exists
because $G$ is locally compact Hausdorff).  Since $\go$ is closed, we may assume $N
\subset \go$ (otherwise replace $N$ by $N \cap \go$).  Since $N \cap
\overline{[z]}$ is compact in $G$ and $N \cap \overline{[z]} = N \cap [z]$, $N
\cap [z]$ is a compact neighborhood of $z$ in $[z]$.

Consider $G{|_{[z]}}=\{\gamma\in G:s(\gamma)\in[z]\text{\ and\ }r(\gamma)\in[z]\}$, the restriction of $G$ to $[z]$.  Since $[z]$
is locally closed, it is locally compact, and hence $G{|_{[z]}}$ is
locally compact.   Thus $G{|_{[z]}}$ satisfies the hypothesis of
\cite[Theorem~2.1]{ramsay}: it is a second countable,  locally compact, Hausdorff
principal groupoid with (one) locally closed orbit $(G{|_{[z]}})^0=[z].$
  Thus $r|:s^{-1}(\{z\}) \rightarrow
[z]$ is a homeomorphism.  (Notice that $s^{-1}(\{z\}) \subseteq G{|_{[z]}}$.) Now
\[
A := r|^{-1}(N \cap [z])=\{\gamma \in G: s(\gamma) =z \text{ and
}r(\gamma) \in N \cap [z]\}
\]
is a compact neighborhood of $z$ in $s^{-1}(\{z\})$.  So $A$ is a compact
subset of $G$ as well. Let $V$ be the interior of $N$ in $\go$. Then
$V$ is an open neighborhood of $z$ in $\go$ such that
\[s^{-1}(\{z\}) \cap r^{-1}(V \cap [z]) = \pi^{-1}((V \cap [z])
\times \{z\})\] is relatively compact in $G$.
\end{proof}

\begin{cor}\label{cor-orbits}
Let $G$ be a second countable, locally compact, Hausdorff, principal groupoid. If
$G$ is integrable then all orbits are locally closed.
\end{cor}

\begin{proof}
We show the contrapositive.      Suppose there exists $z\in G^0$ such that $[z]$ is not
locally closed. Let $V$ be any open relatively compact neighborhood  of $z$ in $G$.  By Lemma~\ref{lem-orbits}
$\lambda_z(r^{-1}(V))=\infty$ and hence \[\sup\{\lambda^x(s^{-1}(\overline{V})):x\in\overline{V}\}\geq \sup\{\lambda_x(r^{-1}(V)):x\in V\}=\infty.\]
So $G$ is not integrable.
\end{proof}

%Note that if $(H,X)$ is an integrable transformation group then the
%orbits are closed; we can improve Corollary~\ref{cor-orbits} in this
%way; the result will follow as a corollary to our main result in the
%sequel.

\begin{prop}\label{prop-int-ktimes}
Let $G$ be a  locally compact, Hausdorff groupoid.
Let $z\in G^0$ and suppose that $G$ fails to be integrable at $z$.
Then there exists a sequence $\{x_n\}$ in $G^0$ such that  $x_n\to z$, and $\{x_n\}$ converges $k$-times in $G^0/G$ to $z$, for every $k\in \P$.  In addition, if  $G$ is second countable, principal and the orbit $[z]$ is locally closed, then $x_n\neq z$ eventually.
\end{prop}
\begin{proof}
Suppose the groupoid fails to be integrable at $z$. Fix $k\in \P$.
Let $\{U_n\}$ be a decreasing sequence of open relatively compact
neighborhoods of $z$ in $G^0$.  By Lemma~\ref{lem-integrable}
\[
\sup_{y\in U_n}\{\lambda^y(s^{-1}(U_n))  \}=\infty
\]
for each $n$. So we can choose a sequence $\{x_n\}$ such that
$x_n\in U_n$ and $\lambda^{x_n}(s^{-1}(U_n))>n$. Note that $x_n\to z$ as $n\to\infty$.

 Let $Q$ be an open symmetric conditionally compact neighborhood of
 $z$ in $G$ and $V$ an  open relatively compact neighborhood of $z$ in $G^0$. By
 Lemma~\ref{lem-measure}(2) there exist $c>0$  such that $\lambda_{v}(Q^2)\leq c$
 whenever $v\in V$.   Choose $n_0$ such that $n_0>(k-1)c$ and $U_{n_0}\subseteq V$.
 Temporarily fix $n>n_0$. Set $\gamma^{(1)}_n=x_n$. For $k\geq 2$ choose
$k-1$ elements $\gamma^{(2)}_n,\dots, \gamma^{(k)}_n$ as
follows.  Note that since $x_n=r(\gamma_n^{(1)})\in V$ we have
\begin{align*}
\lambda_{x_n}\big(r^{-1}(U_n)\setminus Q^2\gamma^{(1)}_n\big)
&\geq \lambda_{x_n}\big(r^{-1}(U_n)\big)-\lambda_{x_n}(Q^2\gamma^{(1)}_n)\\
&=\lambda_{x_n}\big(r^{-1}(U_n)\cap s^{-1}(\{x_n\})\big)-\lambda_{r(\gamma^{(1)}_n)}(Q^2)\\
&> (k-1)c-c=(k-2)c\geq 0.
\end{align*}
So there exists
\[\gamma^{(2)}_n\in \big(r^{-1}(U_n)\cap
s^{-1}(\{x_n\})\big)\setminus Q^2\gamma_n^{(1)};\] note that
$r(\gamma^{(2)}_n)\in U_n\subset V$ and $s(\gamma^{(2)}_n)=x_n$.
Next,
\begin{align*}
\lambda_{x_n}\big((r^{-1}(U_n)\setminus (Q^2\gamma^{(1)}_n\cup
Q^2\gamma^{(2)}_n)\big)
&\geq \lambda_{x_n}(r^{-1}(U_n))-\lambda_{x_n}(Q^2\gamma^{(1)}_n)-\lambda_{x_n}( Q^2\gamma^{(2)}_n)\\
&\geq\lambda_{x_n}(r^{-1}(U_n)\cap s^{-1}(\{x_n\}))-\lambda_{r(\gamma^{(1)}_n)}(Q^2)-\lambda_{r(\gamma^{(2)}_n)}(Q^2)\\
&> (k-3)c\geq 0.
\end{align*}
Continue until $\gamma^{(1)}_n,\dots, \gamma^{(k)}_n$ have been
chosen in this way.

If $n>n_0$, then by construction $s(\gamma^{(i)}_n)=x_n$ and
$r(\gamma^{(i)}_n)\in U_n$ for each $n$; so $r(\gamma^{(i)}_n)\to z$ as $n\to
\infty$  for $1\leq i\leq k$.  Moreover
$\gamma^{(j)}_n(\gamma^{(i)}_n)^{-1}\notin Q^2$ for $1\leq i<j\leq
k$ and $n>n_0$.  To see that
$\{\gamma^{(j)}_n(\gamma^{(i)}_n)^{-1}\}$ tends to infinity, suppose
that it doesn't.  Then,
$\gamma^{(j)}_n(\gamma^{(i)}_n)^{-1}\to\gamma$ by passing to a
subsequence and relabelling.  But then
$s(\gamma^{(j)}_n(\gamma^{(i)}_n)^{-1})=r(\gamma^{(i)}_n)\to z$ and
$r(\gamma^{(j)}_n(\gamma^{(i)}_n)^{-1})=r(\gamma^{(j)}_n)\to z$
implies $\gamma=z$, which is impossible because
$\gamma^{(j)}_n(\gamma^{(i)}_n)^{-1}\notin Q^2$ and $Q$ contains
$G^0$.  Hence $\{x_n\}$ converges $k$-times in $G^0/G$ to $z$.

We claim that if $G$ is second countable and principal, then $x_n\neq z$ eventually. To
see this, suppose $x_n=z$ frequently. Then $\lambda^z(s^{-1}(U_n))>n$ frequently,
and hence
\begin{equation}\label{eq-contra}
\lambda^z(s^{-1}(U_1))=\infty.\end{equation}
The orbits are locally closed and $G$ is second countable and principal, so the source map restricts to a homeomorphism
$s|:r^{-1}(\{z\})\to[z]$.  Since $U_1$ is relatively compact,
$s^{-1}([z] \cap U_1)$
is relatively compact in $r^{-1}(\{z\})$ because  $s|:r^{-1}(\{z\})\to[z]$
is a homeomorphism. But now  $\lambda^z(s^{-1}([z] \cap U_1))=\lambda^z(s^{-1}(U_1))<\infty$, contradicting \eqref{eq-contra}.
\end{proof}

\section{Integrability of $G$ and  trace properties of $C^*(G)$}\label{sec-four}

\begin{prop} \label{prop-converse} Let $G$ be a second-countable, locally compact,
Hausdorff,  principal groupoid.
If $\cs(G)$ has bounded trace then $G$ is integrable.
\end{prop}

The proof  of Proposition~\ref{prop-converse} is based  on that of
\cite[Theorem~2.3]{MW90}.  There, Muhly and Williams choose a sequence $\{x_n\}\subseteq G^0$
with $x_n\to z$ which witnesses the failure of the groupoid to be proper.
They then carefully construct a function $f\in C_c(G)$ to obtain an element $d$
of the Pedersen ideal of $\cs(G)$ such that $\tr(L^{x_n}(d))$ does not converge to
$\tr(L^z(d))$. Since the Pedersen ideal is the minimal dense ideal \cite[Theorem~5.6.1]{ped-auto}, the ideal of continuous-trace elements cannot be dense, so $\cs(G)$ does not have continuous trace. We adopt the same strategy, use exactly the same function
$f$, but adapt the proof of \cite[Theorem~2.3]{MW90} using ideas from \cite[Proposition~3.5]{aH02}.

\begin{proof}[Proof of Proposition~\ref{prop-converse}]
Fix $M\in\P$. We will show that there is an element $d$ of the Pedersen
ideal of $C^*(G)$, a sequence of representations $\{L^{x_n}\}$ and $n_0>0$ such that
$\tr(L^{x_n}(d))>M$ whenever $n>n_0$.  Since $M$ is arbitrary, $C^*(G)$ cannot have bounded trace.

If $G$ is not integrable, then the integrability fails at some $z\in G^0$ by
Lemma~\ref{lem-integrable}.  If the orbits are not closed then $C^*(G)$ cannot
be CCR by \cite[Theorem~4.1]{C}, hence cannot have bounded trace.  So from now on we may assume that
the orbits are closed. By Proposition~\ref{prop-int-ktimes}, there exists a sequence $\{x_n\}$
such that $x_n\neq z$, $x_n\to z$, and $\{x_n\}$ converges $k$-times in $G^0/G$ to $z$, for every $k\in\P$.

Since we will use exactly the same function $f$ that was used in the proof of \cite[Theorem~2.3]{MW90}, our first task is to briefly outline its construction. Fix a function $g \in \cc(G^0)$
such that $0\leq g\leq 1$ and $g$ is identically one on a neighborhood $U$ of $z$.
Let $N = \supp  g$ and
\begin{align*}
F_z^N &:= s^{-1}(\{z\}) \cap r^{-1}([z] \cap N)=  s^{-1}(\{z\}) \cap r^{-1}( N)\\
F^z_N &:= r^{-1}(\{z\}) \cap s^{-1}([z] \cap N)= r^{-1}(\{z\}) \cap s^{-1}(N).
\end{align*}
There exist symmetric, open, conditionally compact neighborhoods $W_0$ and $W_1$ in $G$ such that
\[
G^0 \subseteq W_0 \subseteq \overline{W_0} \subseteq W_1\quad\text{and}\quad F^z_N \cup F^N_z\subseteq W_0.
\]
Thus
$
\overline{W}_1^7 z \setminus W_0 z \subseteq r^{-1}(G^0 \setminus N)
$.
(The reason for using $\overline{W_1}^7$ becomes clear at \eqref{eq-w1to7} below).
By a compactness argument, there exist open, symmetric, relatively compact
neighborhoods $V_0\subseteq G^0$ and $V_1$ of $z$ in $G$ such that $\overline{V_0}\subset V_1$ and
\begin{equation}
\overline{W}^7_1 \overline{V}_0 \setminus W_0 V_0 \subseteq r^{-1}(G^0 \setminus N).
\label{bill}
\end{equation}
Now note that if $\gamma\in\overline{W}^7_1 \overline{V}_1 \overline{W}^7_1\setminus W_0V_0W_0$ then $r(\gamma)\in r(\overline{W}^7_1 \overline{V}_0 \setminus W_0 V_0)\subseteq G^0\setminus N$. It  follows that the function $g^{(1)}:G\to [0,1]$ defined by
\begin{equation}
g^{(1)}(\gamma) =
\begin{cases}
g(r(\gamma))& \text{if }\gamma \in \overline{W}^7_1 \overline{V}_1 \overline{W}^7_1,\\
0&\text{if } \gamma \notin W_0V_0W_0
\end{cases}
\notag
\end{equation}
is well-defined and continuous with compact support in $G$.
By construction
\begin{equation}
(W_0V_0W_0)^2 = W_0V_0W_0^2V_0W_0 \subseteq W_0^4V_0W_0^4 \subseteq
\overline{W}_0^4\overline{V}_0\overline{W}_0^4 \subseteq W_1^4V_1W_1^4 \subseteq
\overline{W}^4_1\overline{V}_1\overline{W}^4_1.
\notag
\end{equation}
So there exists a function $b\in \cc(G)$ such that $0 \leq b \leq 1$, $b$ is
identically one on $W_0V_0W_0^2V_0W_0$ and is identically zero on the complement
of $\overline{W}^4_1\overline{V}_1\overline{W}^4_1.$  Further, we can replace $b$
with $(b+b^*)/2$ to ensure that $b$ is self-adjoint. Set
\begin{equation*}
f(\gamma) = g(r(\gamma))g(s(\gamma))b(\gamma);
\end{equation*}
note that $f \in \cc(G)$  is self-adjoint.

For $\xi \in L^2(G,\lambda_{u})$ and $\gamma\in G$ we have
\begin{align}
L^{u}(f)\xi(\gamma) &= \int  f(\gamma \alpha)\xi(\alpha^{-1})\ d\lambda^{u}(\alpha)\notag\\
&=\int  g(r(\gamma))g(s(\alpha))b(\gamma\alpha)\xi(\alpha^{-1})\ d\lambda^{u}(\alpha)\notag\\
&=g(r(\gamma))\int  g(s(\alpha))b(\gamma\alpha)\xi(\alpha^{-1})\ d\lambda^{u}(\alpha)\notag\\
&=g(r(\gamma))\int  g(r(\alpha))b(\gamma\alpha^{-1})\xi(\alpha)\ d\lambda_{u}(\alpha)\label{rep-formula}.
\end{align}
By \cite[Lemma~2.8]{MW90},  $g^{(1)}$ is an eigenvector for $L^{x_n}(f)$ with eigenvalue
\[
\mu_{x_n}^{(1)}=\int  g(r(\alpha))g^{(1)}(\alpha) \ d\lambda_{x_n}(\alpha)
=\int_{W_0V_0W_0}g(r(\alpha))^2 \ d\lambda_{x_n}(\alpha).
\]
By \cite[Lemma~2.9]{MW90}, there exist an open $V_2 \subseteq V_0$ and a conditionally compact neighborhood
$Y$ of $G^0$ so that $Y \subseteq W_0$ and if $v \in V_2$, then $r(Yv) \subseteq U$.  Notice that $YV_2Y $
is a relatively compact subset of $W_0V_0W_0$.
By Lemma~\ref{lem-measure}(1) there exist an open  neighborhood $V_3$ of $z$ and $a>0$ such that
\begin{equation}\label{eq-a}\lambda_v(YV_2Y)\geq a \text{\ whenever\ } v\in V_3.\end{equation}
Now, if $\alpha\in YV_2Y$ then $r(\alpha)\in U$ and hence $g(r(\alpha))=1$; it follows that
\begin{equation*}
\mu_{x_n}^{(1)} \geq \int_{YV_2Y} g(r(\alpha))^2 \ d\lambda_{x_n}(\alpha)=\lambda_{x_n}(YV_2Y)\geq a>0
\end{equation*}
whenever $x_n\in V_3$.

So far our set-up is the one from \cite{MW90}.  Now choose $l\in\P$ such that $la^2>M$.
(Note that $a$ is independent of $l$!) The sequence $\{x_n\}$ converges $k$-times in $G/G^0$ to
$z$ for every $k\in\P$, so it certainly converges $l$ times. So
there exist $l$ sequences
\[
\{ \gamma_n^{(1)}\},\  \{\gamma_n^{(2)}\},\dots,\{\gamma_n^{(l)}\}\subseteq G
\]
such that
\begin{enumerate}
\item $r(\gamma_n^{(i)})\to z$ as $n\to\infty$ for $1\leq i\leq l$;
\item  $s(\gamma_n^{(i)})=x_n$ for $1\leq i\leq k$;
\item if $1\leq i<j\leq l$ then $\gamma_n^{(j)}(\gamma_n^{(i)})^{-1}\to\infty$.
\end{enumerate}
Moreover, by construction (see Proposition~\ref{prop-int-ktimes}), we may
take $\gamma_n^{(1)}=x_n$. Temporarily fix $n$.  Set
$g_n^{(1)}:=g^{(1)}$, and for $2 \leq j \leq l$ set
\begin{align*}
g_n^{(j)}(\gamma) &:=
\begin{cases}
g^{(1)}(\gamma(\gamma^j_n)^{-1}),&\text{if }s(\gamma)=s(\gamma^j_n);\\
0,& \text{otherwise}
\end{cases}\\
&=
\begin{cases}
g(r(\gamma)),& \text{if } \gamma \in \overline{W}^7_1 \overline{V}_1 \overline{W}^7_1\gamma_n^j;\\
0,&\text{otherwise}
\end{cases}\\
&=
\begin{cases}
g(r(\gamma)),& \text{if }\gamma \in \overline{W}^7_1 \overline{V}_1 \overline{W}^7_1\gamma_n^j;\\
0, &\text{if } \gamma \notin W_0V_0W_0\gamma^j_n.
\end{cases}
\end{align*}
Each $g_n^{(i)}$ $(1\leq j\leq l)$ is a well-defined function in
$C_c(G)$ with support contained in $W_0V_0W_o\gamma_n^{(j)}$. For $1\leq i<j\leq l$,  $\gamma_n^{(j)}(\gamma_n^{(i)})^{-1}\notin (W_0V_0W_0)^2$ eventually, so there exists $n_0>0$ such that, for every $0\leq i,j\leq l$, $i\neq j$,
\[W_0V_0W_0\gamma_n^{(j)}\cap W_0V_0W_0\gamma_n^{(i)}=\emptyset\]
whenever $n>n_0$.

We now prove a generalization of \cite[Lemma~2.8]{MW90}
which, together with \eqref{rep-formula}, immediately implies that
each $g_n^{(j)}$ is an eigenvector of $L^{x_n}(f)$ for $1\leq j\leq
l$.

\begin{lemma}\label{new2.8}
With the choices made above, for all $\alpha, \gamma \in G$ and $1\leq j\leq l$,
\begin{equation*}
g(r(\gamma))g(r(\alpha))b(\gamma\alpha^{-1})g_n^{(j)}(\alpha) = g_n^{(j)}(\gamma)g(r(\alpha))g_n^{(j)}(\alpha).
\end{equation*}
\end{lemma}
\begin{proof}
If $\alpha \notin W_0V_0W_0\gamma^{(j)}_n$, then both sides are zero.  So we may assume
throughout that $\alpha \in W_0V_0W_0\gamma^{(j)}_n$.

If $\gamma \in W_0V_0W_0\gamma^{(j)}_n$, then $g_n^{(j)}(\gamma)=g(r(\gamma))$ and
$\gamma\alpha^{-1} \in W_0V_0W_0^2V_0W_0$, so $b(\gamma\alpha^{-1})=1$ and both sides agree.

If $\gamma \in \overline{W}^7_1\overline{V}_1\overline{W}_1^7\gamma^{(j)}_n \setminus W_0V_0W_0\gamma_n^{(j)}$, then
$g(r(\gamma))=0=g_n^{(j)}(\gamma)$, so both sides are zero.

Finally, if $\gamma \notin \overline{W}^7_1\overline{V}_1\overline{W}_1^7\gamma^{(j)}_n$,
then $g_n^{(j)}(\gamma)=0$ so the right-hand side zero.  On the other hand,
if $\gamma\alpha^{-1} \in\overline{W}_1^4\overline{V}_1\overline{W}_1^4(=\supp b)$ then
\begin{equation}\gamma \in\overline{W}_1^4\overline{V}_1\overline{W}_1^7 \gamma_n^{(j)}
\subseteq \overline{W}_1^7\overline{V}_1\overline{W}_1^7\gamma^{(j)}_n.\label{eq-w1to7}
\end{equation}
So $\gamma \notin \overline{W}^7_1\overline{V}_1\overline{W}_1^7\gamma^{(j)}_n$ implies $\gamma\alpha^{-1} \notin \supp b$, so the left-hand side is zero as well.
\end{proof}
Let $\mu_n^{(j)}$  be the eigenvalue  corresponding to the eigenvector $g_n^{(j)}$.
Using \eqref{eq-a}
\[
\mu_n^{(j)} =\int _{W_0V_0W_0\gamma^{(j)}_n } g(r(\alpha))^2 \,d\lambda_{x_n}(\alpha)
\geq \lambda_{x_n}(YV_2Y\gamma^{(j)}_n)=\lambda_{r(\gamma^{(j)}_n)}(YV_2Y)\geq a
\]
whenever $r(\gamma^{(j)}_n)\in V_3$.
Choose $n_1>n_0$ such that $n>n_1$ implies $x_n\in V_3$ and $r(\gamma^{(j)}_n)\in V_3$ for $1\leq j\leq l$. Then $L^{x_n}(f*f)$ is a positive compact operator with  $l$ eigenvalues $\mu_n^{(j)}\geq a^2$ for $1\leq j\leq l$. To push $f*f$ into the Pedersen ideal, let $r\in C_c(0,\infty)$ be any function satisfying
\[
r(t)=\begin{cases}
0,& \text{if }t<\frac{a^2}{3};\\
2t-\frac{2a^2}{3}, &\text{if }\frac{a^2}{3}\leq t< \frac{2a^2}{3};\\
t,&\text{if }\frac{2a^2}{3}\leq t\leq \|f*f\|.
\end{cases}
\]
Set $d:=r(f*f)$.  Now $d$ is a positive element of
the Pedersen ideal of $C^*(G)$ with $\tr(L^{x_n}(d))\geq la^2>M$ whenever $n>n_1$. Since $M$ was arbitrary, $L^x\mapsto \tr(L^{x}(d))$ is unbounded on $C^*(G)^\wedge$.
Thus $\cs(G)$ does not have bounded trace.
\end{proof}

\begin{prop}
\label{prop-trace} Suppose $G$ is a second
countable, locally compact, Hausdorff,  principal groupoid.  If $G$ is integrable then $\cs(G)$
has bounded trace.
\end{prop}

\begin{proof}
Since $G$ is principal and integrable,   the
orbits are locally closed by Corollary~\ref{cor-orbits}, and $x\mapsto L^x$ induces a homeomorphism of  $\go/G$ onto $\cs(G)^{\wedge}$ by
\cite[Proposition~5.1]{C}.  To  show that $\cs(G)$ has bounded trace,
it suffices to see that for a fixed $u\in G^0$ and all $f \in C_c(G)$,  $\tr(L^u(f^**f))$ is
bounded independent of $u$.

Fix $u\in G^0$ and let $\xi \in L^2(G,\lambda_u)$. Since
\[
L^u(f)\xi(\gamma) = \int  f(\gamma \alpha^{-1}) \xi(\alpha)\
d\lambda_u(\alpha),
\]
$L^u(f)$ is a kernel operator on $L^2(G,\lambda_u)$ with kernel
$k_f$ given by
$
k_f (\gamma, \alpha) = f(\gamma \alpha^{-1})
$.
We will show that $k_f \in L^2(G \times G, \lambda_u \times \lambda_u)$ and find a bound on $k_f$ independent of $u$.  This will complete the proof  since
$
\tr(L^u(f^**f))=\|k_f\|^2
$
by, for example, \cite[Theorem~3.4.16]{pedersen}.

Notice that
\begin{align}
\|k_f\|^2=\int_{G \times G} |k_f(\gamma, \alpha)|^2 \ d(\lambda_u\times\lambda_u)(\gamma,\alpha)
&= \int_G   \int_G  |f(\gamma \alpha^{-1})|^2  \
d\lambda_u(\gamma) \ d\lambda_u(\alpha)\notag \\
&=\int_G \int_G  |f(\gamma)|^2 \
d\lambda_{r(\alpha)} (\gamma) \ d\lambda_u(\alpha) \label{trace}
\end{align}
by Tonelli's Theorem and right invariance.
For a fixed $\alpha$, the inner integral
\[
\int_G  |f(\gamma)|^2 \ d\lambda_{r(\alpha)}(\gamma) \leq \|f\|_{\infty}^2
\lambda_{r(\alpha)}(\supp f)
\]
and is zero unless $r(\alpha) \in s(\supp f)$.  The outer integral
is zero unless $s(\alpha)=u$.  Let
\[
K =r^{-1}(s(\supp f))\cap s^{-1}(\{u\}).
\]
So
\begin{align*}
\eqref{trace}
&\leq \int_K \|f\|_{\infty}^2 \lambda_{r(\alpha)}(\supp f) \
d\lambda_u(\alpha)\\
&\leq \|f\|_{\infty}^2  \sup
\big\{\lambda_{r(\alpha)}(\supp f):r(\alpha) \in s(\supp f)\big\} \lambda_u(K)\\
&\leq \|f\|_{\infty}^2  \sup
\big\{\lambda_x(\supp f):x \in s(\supp f)\big\} \sup\big\{\lambda_x(r^{-1}(s(\supp f)):x\in G^0\big\}.
\end{align*}

Since $s(\supp f)$ is a compact subset of $G^0$, by integrability there exists $N>0$ such that
\[
\sup\big\{\lambda_x(r^{-1}(s(\supp f)):x\in G^0\big\}<N
\]
(see also Remark~\ref{rem-integr}). 
Note that $N$ does not depend on $u$.
By Lemma~\ref{lem-measure}(2), applied to the conditionally compact neighborhood $\supp f$ and the relatively compact neighborhood $s(\supp f)$, there exists $M>0$ such that $\lambda_x(\supp f)<M$ for all $x\in s(\supp f)$, that is,
\[\sup\{\lambda_x(\supp f):x\in s(\supp f)\}<M.\]
Note that $M$ does not depend on $u$.

Thus $\|k_f\|^2<\|f\|^2_\infty MN$, so
$k_f \in L^2(G\times G, \lambda_u \times \lambda_u)$ as claimed, and 
\[
\tr(L^u(f^**f))=\|k_f\|^2< \|f\|^2_\infty M N,
\]
which is a bound on $\tr(L^u(f^**f))$ independent of $u$.
\end{proof}

Combining Propositions~\ref{prop-converse} and~\ref{prop-trace} we have:

\begin{thm}\label{thm-main}
Suppose $G$ is a second countable, locally compact, Hausdorff,
principal groupoid.  Then $G$ is integrable if and only if
$\cs(G)$ has bounded trace.
\end{thm}

By Corollary~\ref{cor-orbits}, integrable groupoids have locally closed orbits; we can now improve this result using Proposition~\ref{prop-trace} (which relied on  Corollary~\ref{cor-orbits}).  The $C^*$-algebras of integrable groupoids have bounded trace by Propositon~\ref{prop-trace}, hence are CCR. But if $C^*(G)$ is CCR then the orbit space $G^0/G$ is $T_1$ by \cite[Theorem~4.1]{C}, that is, the orbits are closed.  

\begin{cor}
Suppose $G$ is a second countable, locally compact, Hausdorff,
principal groupoid.  If $G$ is integrable then the orbits of
$G$ are closed.
\end{cor}
We have so far been unable to prove directly that integrability implies that the orbits are closed.


\begin{thebibliography}{99}

\bibitem{AD} R.J. Archbold and K. Deicke, \emph{Bounded trace $C^*$-algebras and integrable actions},
Math. Zeit. \textbf{250} (2005), 393--410.

\bibitem{AaH} R.J. Archbold and A. an Huef, \emph{Strength of convergence in the orbit space of a transformation group},
to appear in J. Funct. Anal.

\bibitem{AS}  R.J. Archbold and D.W.B. Somerset, \emph{Transition probabilities and trace functions for $C\sp *$-algebras},  Math. Scand.  \textbf{73}  (1993),  81--111. 

\bibitem{C} L.O. Clark, \emph{Classifying the type of principal groupoid $C^*$-algebras}, to appear in J. Operator Theory.

\bibitem{green}
P. Green, \emph{{$C^*$}-algebras of transformation groups with smooth orbit
  space}, Pacific J. Math. \textbf{72} (1977), 71--97.

\bibitem{aH02}
A.~an Huef, \emph{Integrable actions and the transformation groups whose
$C^*$-algebras have bounded trace}, Indiana Univ. Math. J. \textbf{51} (2002), 1197--1233.

\bibitem{MRW87}
P.~S. Muhly, J. Renault, and D.~P. Williams, \emph{Equivalence and
  isomorphism for groupoid {\cs}-algebras}, J. Operator Theory
\textbf{17}
  (1987), 3--22.

\bibitem{MW90} P.S. Muhly and D.P. Williams, \emph{Continuous trace groupoid $C^*$-algebras}, Math. Scand. \textbf{66} (1990), 231--241.

\bibitem{ped-auto}
G.K. Pedersen, \emph{$C^*$}-algebras and their automorphism groups, Academic
  Press, London, 1979.

\bibitem{pedersen} G.~K. Pedersen, \emph{Analysis Now},
Springer-Verlag, New York, 1989.

\bibitem{ramsay}
A. Ramsay, \emph{The Mackey-Glimm dichotomy for foliations and other
Polish groupoids}, J. Funct. Anal. \textbf{94} (1990), 358--374.

\bibitem{renault}
J. Renault, \emph{ A groupoid approach to $C^*$-algebras}, Lecture
Notes in Mathematics, No. \textbf{793}, Springer-Verlag, New York,
1980.

\bibitem{rieffel} M.A. Rieffel, \emph{Integrable and proper actions on $C^*$-algebras, and square-integrable representations of groups}, Expo. Math. \textbf{22} (2004), 1--53.


\end{thebibliography}
\end{document}